\pdfoutput=1
\RequirePackage{ifpdf}
\ifpdf 
\documentclass[pdftex]{sigma}
\else
\documentclass{sigma}
\fi

\newcommand{\qbinom}[2]{\genfrac[]{0pt}0{#1}{#2}}

\newcommand{\boi}{\mathbf i}
\newcommand{\boj}{\mathbf j}
\newcommand{\bop}{\mathbf p}
\newcommand{\bc}{\mathbf C}
\newcommand{\bn}{\mathbf N}
\newcommand{\bu}{\mathbf U}
\newcommand{\bz}{\mathbf Z}

\newcommand{\ev}{\operatorname{ev}}
\newcommand{\wt}{\operatorname{wt}}
\newcommand{\ch}{\operatorname{ch}}
\newcommand{\gr}{\operatorname{gr}}
\newcommand{\loc}{\operatorname{loc}}
\newcommand{\Hom}{\operatorname{Hom}}
\newcommand{\res}{\operatorname{res}}

\numberwithin{equation}{section}

\newtheorem{Theorem}{Theorem}[section]
\newtheorem{Lemma}[Theorem]{Lemma}
\newtheorem{Proposition}[Theorem]{Proposition}
{\theoremstyle{definition}

\newtheorem{Remark}[Theorem]{Remark}
}

\begin{document}

\allowdisplaybreaks

\renewcommand{\thefootnote}{$\star$}

\renewcommand{\PaperNumber}{032}

\FirstPageHeading

\ShortArticleName{Modules with Demazure Flags and Character Formulae}

\ArticleName{Modules with Demazure Flags\\ and Character Formulae\footnote{This paper is a~contribution to the
Special Issue in honor of Anatol Kirillov and Tetsuji Miwa.
The full collection is available at \href{http://www.emis.de/journals/SIGMA/InfiniteAnalysis2013.html}
{http://www.emis.de/journals/SIGMA/InfiniteAnalysis2013.html}}}

\Author{Vyjayanthi CHARI, Lisa SCHNEIDER, Peri SHEREEN and Jeffrey WAND}

\AuthorNameForHeading{V.~Chari, L.~Schneider, P.~Shereen and J.~Wand}

\Address{Department of Mathematics, University of California, Riverside, CA 92521, USA}
\Email{\href{mailto:chari@math.ucr.edu}{chari@math.ucr.edu}, \href{mailto:lschn005@ucr.edu}{lschn005@ucr.edu}, \href{mailto:psher001@ucr.edu}{psher001@ucr.edu}, \href{mailto:wand@math.ucr.edu}{wand@math.ucr.edu}}

\ArticleDates{Received October 22, 2013, in f\/inal form March 17, 2014; Published online March 27, 2014}

\Abstract{In this paper we study a~family of f\/inite-dimensional graded representations of the current algebra of $\mathfrak{sl}_2$ which are indexed by partitions.
 We show that these representations admit a~f\/lag where the successive quotients are Demazure modules which occur in a~level $\ell$-integrable module for $A_1^1$ as long as $\ell$ is large. We associate to each partition and to each $\ell$ an edge-labeled directed graph which allows us to describe in a~combinatorial way the graded multiplicity of a~given level $\ell$-Demazure module in the f\/iltration. In the special case of the partition $1^s$ and $\ell=2$, we give a~closed formula for the graded multiplicity of level two Demazure modules in a~level one Demazure module. As an application, we use our result along with the results of Naoi and Lenart et al., to give the character of a~$\mathfrak{g}$-stable level one Demazure module associated to $B_n^1$ as an explicit combination of suitably specialized Macdonald polynomials. In the case of $\mathfrak{sl}_2$, we also study the f\/iltration of the level two Demazure module by level three Demazure modules and compute the numerical f\/iltration multiplicities and show that the graded multiplicites are related to (variants of) partial theta series.}


\Classification{06B15 ; 05E10; 14H42}

\renewcommand{\thefootnote}{\arabic{footnote}}
\setcounter{footnote}{0}

\vspace{-2mm}

\section{Introduction}

The current algebra associated to a~simple Lie algebra $\mathfrak{g}$ is just the Lie algebra of polynomial maps $\bc\to \mathfrak{g}$ or equivalently, the special maximal parabolic subalgebra of the untwisted af\/f\/ine Lie algebra $\hat{\mathfrak{g}}$ associated to $\mathfrak{g}$. Given any integrable highest weight representation of the af\/f\/ine Lie algebra, one can def\/ine a~certain family of modules for the appropriate Borel subalgebra in~$\hat{\mathfrak{g}}$. If the center of the af\/f\/ine algebra acts by the positive integer~$\ell$, these modules are referred to as the Demazure module of level~$\ell$. In certain special cases, the Demazure modules are actually modules for the parabolic subalgebra, i.e., for the current algebra and from now on, we shall only be concerned with such Demazure modules.

The notion of a~local Weyl module for the current algebra was introduced and f\/irst studied in~\cite{CPweyl} and the def\/inition was motivated by the representation theory of quantum af\/f\/ine algebras. It was proved later (see \cite{CL, CV,FoL}) that in the simply-laced case a~local Weyl module is isomorphic to a~Demazure module. In the non-simply laced case, it was proved in \cite{Naoi} that a~local Weyl module has a~f\/iltration by Demazure modules. One of the important ingredients in \cite{Naoi} is the following: a~level $\ell$-Demazure module for~$A_n^1$ has a~f\/iltration by level $m$ Demazure modules for all $m\ge\ell$. The proof of this statement requires taking the $q=1$ limit of
a~result due to A.~Joseph~\cite{Jo} in the quantum case using the theory of canonical bases. However, Joseph's result was only proved for simply-laced Kac--Moody Lie algebras and hence did not include $A_1^1$. It is remarked in \cite[Remark~4.15]{Naoi} that an inspection of Joseph's proof along with a~positivity result proved in \cite{L} showed that Joseph's proof worked worked for $A_1^1$ as well. Our paper was motivated by an attempt to see if it was possible to give a~direct proof of this result for $A_1^1$. This is achieved as a~particular case of Theorem~\ref{main1}. Our methods also allow us to determine the multiplicity of the level two Demazure module in the local Weyl module. Using \cite[Sections~4.19 and~9]{Naoi} this means that we give the f\/iltration multiplicities (see Proposition~\ref{chdem}) for a~Demazure module occurring in the f\/iltration of a~local Weyl module for~$B_n$. In Subsection~3.10, 
we study the case of the Demazure module of level two and show that the multiplicities of the level three Demazure f\/lag are related to partial theta series.

The paper is organized as follows. In Section~\ref{section2}, we recall the def\/initions and collect the necessary results from \cite{CL,FoL,Naoi} on the relationship between the $\mathfrak{g}$-stable level one Demazure modules and local Weyl modules. In the next section we restrict our attention to~$\mathfrak{sl}_2$ and recall the def\/inition from~\cite{CV} of the modules~$V(\xi)$. We prove that these modules admit a~Demazure f\/lag of level $\ell$ for all $\ell$ suf\/f\/iciently large. In the case when $\ell=2$ we give a~closed formula for the multiplicities. We then use results in~\cite{Naoi} to give a~character formula for the local Weyl module in terms of the Demazure module when~$\mathfrak{g}$ is of type~$B_n$. Together with the results of~\cite{LNSSS} we can give the character of the Demazure module for~$B_n^1$ as a~combination of Macdonald polynomials. In the last section, we def\/ine a~directed, edge labeled, graph $\mathbb H_\ell(\xi)$ and show that the multiplicity of a~level~$\ell$ Demazure module in~$V(\xi)$ is given by natural polynomials associated to specif\/ic paths of this graph.

\section[Local Weyl modules and Demazure modules of level $\ell$]{Local Weyl modules and Demazure modules of level $\boldsymbol{\ell}$}\label{section2}

We summarize the main def\/initions and results needed for our paper. We refer the reader to \cite{CFK} for further information on the local Weyl modules, to \cite{CL,FoL,Naoi} for the def\/inition of Demazure modules and their connection with local Weyl modules.

{\bf 2.1.}~Throughout this paper $\bc$ denotes the f\/ield of complex numbers and $\bz$, $\bz_+$, $\bn$ denote the integers, non-negative integers and positive integers respectively. Given a~complex Lie algebra~$\mathfrak{a}$ denote by~$\bu(\mathfrak{a})$ the corresponding universal enveloping algebra. The associated current Lie algebra is denoted $\mathfrak{a}[t]$: as a~vector space it is just $\mathfrak{a}\otimes \bc[t]$ and the Lie bracket is given in the natural way (i.e., $[a\otimes f, b\otimes g]=[a,b]\otimes f g$, for all $a,b\in\mathfrak{a}$ and $f,g\in\bc[t]$). (Here $\bc[t]$ is the polynomial algebra in an indeterminate $t$). The degree grading on $\bc[t]$ def\/ines a~natural $\bz_+$-grading on $\mathfrak{a}[t]$ and hence also on $\bu(\mathfrak{a}[t])$: an element of the form $(a_1\otimes t^{r_1})\cdots (a_s\otimes t^{r_s})$ has grade $r_1+\cdots+r_s$.

 A {\em graded} representation of $\mathfrak{a}[t]$ is a~$\bz$-graded vector space which admits a~compatible Lie algebra action of $\mathfrak{a}[t]$, i.e.,
\[
V=\bigoplus_{r\in\bz}V[r],\qquad (\mathfrak{a}\otimes t^s)V[r]\subset V[r+s],\qquad r\in\bz, \qquad s\in\bz_+.
\] If $V$ and $V'$ are graded $\mathfrak{a}[t]$-modules, we say that $\pi:V\to V'$ is a~morphism of graded $\mathfrak{a}[t]$-modules if $\pi$ is a~degree zero morphism of $\mathfrak{a}[t]$-modules. In particular, this means that $V[r]$ is a~$\mathfrak{a}$-submodule of $V$. For $r\in\bz$, let $\tau_r$ be the grading shift operator: if $V$ is a~graded $\mathfrak{a}[t]$-module, then $\tau_rV$ is the graded $\mathfrak{a}[t]$-module with the graded pieces shifted uniformly by $r$ and the action of $\mathfrak{a}[t]$ unchanged. If $M$ is an $\mathfrak{a}$-module
 def\/ine a~$\mathfrak{a}[t]$-module structure on $M$ by: $(a\otimes t^r)m= \delta_{r,0}am$ where $a\in\mathfrak{a}$ and $r\in\bz_+$ and $\delta_{r,0}$ is the Kronecker delta symbol. We denote this module as $\ev_0 M$ and observe that $\ev_0M$ is a~graded $\mathfrak{a}[t]$-module with $\left(\ev_0M\right)[0]=M$.

{\bf 2.2.}~From now on, $\mathfrak{g}$ is an arbitrary simple f\/inite-dimensional complex Lie algebra of rank~$n$, $\mathfrak{h}$~is a~f\/ixed Cartan subalgebra of $\mathfrak{g}$ and $R$~is the set of roots of $\mathfrak{g}$ with respect to $\mathfrak{h}$. The restriction of the Killing form of $\mathfrak{g}$ to $\mathfrak{h}$ induces an isomorphism between $\mathfrak{h}$ and $\mathfrak{h}^*$ and hence also a~symmetric non-degenerate form $(\ , \ )$ on $\mathfrak{h}^*$. We shall assume that this form on $\mathfrak{h}^*$ is normalized so that the square length of a~long root is two and for $\alpha\in R$ set $d_\alpha=\frac{2}{(\alpha,\alpha)}$.

 Let $I=\{1,\dots ,n\}$ and f\/ix a~set $\{\alpha_i: i\in I\}$ of simple roots for $R$ and let $\{\omega_i: i\in I\}\subset\mathfrak{h}^*$ be the set of fundamental weights. Let $Q$ (resp.~$Q^+$) be the $\bz$ span (resp.\ the $\bz_+$ span) of $\{\alpha_i: i\in I\}$ and similarly def\/ine $P$ (resp.\ $P^+$) to be the $\bz$ (resp.~$\bz_+$) span of $\{\omega_i: i\in I\}$ and set $R^+=R\cap Q^+$. Def\/ine a~partial order on $P^+$ by $\lambda\ge\mu$ if and only if $\lambda-\mu\in Q^+$.

 Finally, let $\{x_\alpha^\pm, h_i: \alpha\in R^+, i\in I\}$ be a~Chevalley basis for $\mathfrak{g}$ and set $x_i^\pm=x_{\alpha_i}^\pm$ for $i\in I$.

{\bf 2.3.}~For $\mu\in P^+$, let $V(\mu)$ be the irreducible f\/inite-dimensional $\mathfrak{g}$-module generated by an element $v_\mu$ with def\/ining relations
\[
x_i^+v_\mu=0,\qquad h_iv_\mu=\mu(h_i)v_\mu,\qquad(x_i^-)^{\mu(h_i)+1}v_\mu=0,\qquad i\in I.
\]
It is well-known that any f\/inite-dimensional $\mathfrak{g}$-module $V$ is isomorphic to a~direct sum of irreducible modules $V(\mu)$, $\mu\in P^+$. Further, we may write
\[
V=\bigoplus_{\nu\in P} V_\nu,\qquad V_\nu=\{v\in V: hv=\nu(h)v,\; h\in\mathfrak{h}\},
\]
and we set $\wt V=\{\nu\in P: V_\nu\ne 0\}$. Let $\bz[P]$
denote the group ring of $P$ with integer coef\/f\/icients and basis $e(\mu)$, $\mu\in P$. The character of a~f\/inite-dimensional $\mathfrak{g}$-module $V$ is the element of $\bz[P]$ def\/ined by
\[
\ch_{\mathfrak{g}}(V)=\sum_{\mu\in P} \dim V_{\mu}e(\mu).
\] Since $\wt V(\mu)\subset \mu-Q^+$, it follows that the
characters $\ch_{\mathfrak{g}}V(\mu)$, $\mu\in P^+$ form a~linearly independent subset of $\bz[P]$.
Given a~graded $\mathfrak{g}[t]$-module and an indeterminate $q$, we let
\begin{gather*}
\ch_{\gr}(V)=\sum_{r\ge 0}\ch_{\mathfrak{g}}V[r]q^r\in\bz[P][q].
\end{gather*}

{\bf 2.4.}~For $\lambda\in P^+$, the local Weyl module $W_{\loc}(\lambda)$ is the $\mathfrak{g}[t]$-module generated by an element~$w_\lambda$ with def\/ining relations, for $i\in I$ and $s\in\bz_+$,
\begin{gather*}
(x^+_i\otimes\bc[t])w_\lambda=0,\qquad (h_i\otimes t^s)w_\lambda=\lambda(h_i)\delta_{s,0}w_{\lambda},\qquad
(x_i^-\otimes 1)^{\lambda(h_i)+1}w_\lambda=0.
\end{gather*}
It is trivial to see that $\wt W_{\loc}(\lambda)\subset \lambda-Q^+$.
 The local Weyl module is clearly $\bz_+$-graded once we declare the grade of $w_\lambda$ to be zero, and
\begin{gather*}W_{\loc}(\lambda)[0]\cong_{\mathfrak{g}} V(\lambda),
\qquad \ch_{\gr} W_{\loc}(\lambda)=\ch_{\mathfrak{g}} V(\lambda)+ \!\!\sum_{\stackrel{(\mu,r)\in P^+\times \bn}{\mu< \lambda}}\!\!\!
 \dim\Hom_{\mathfrak{g}}(V(\mu), W_{\loc}(\lambda)[r])q^r.
 \end{gather*}
Finally, $\ev_0 V(\lambda)$ is the unique graded irreducible quotient of $W_{\loc}(\lambda)$
and any irreducible graded f\/inite-dimensional module is isomorphic to $\tau_s\ev_0 V(\lambda)$ for some $(\lambda,s)\in P^+\times\bz_+$.

{\bf 2.5.}~The Demazure module of level $\ell$ and weight $\lambda$ is denoted by $D(\ell,\lambda)$ and for our purposes the def\/inition is as follows: it is the graded quotient of $W_{\loc}(\lambda)$ by the submodule generated by the elements
 \begin{gather} \label{demrelo}
 \{(x^-_\alpha\otimes t^p)^{r+1}w_\lambda:\ p\in\bz_+,\ r\ge\max\{0,\lambda(h_\alpha)-d_\alpha\ell p\},\ \text{for} \ \alpha\in R^+\}.
 \end{gather} By abuse of notation, we continue to denote the image of $w_\lambda$ in $D(\ell,\lambda)$ by $w_\lambda$.
 Notice that for all $\alpha\in R^+$, $p\geq 1$, $(x^-_{\alpha}\otimes t^p)w_\lambda=0$ if $\ell$ is such that $\ell\ge\lambda(h_\alpha)$.

 Again, we have
 \begin{gather}\label{demeva} D(\ell,\lambda)[0]\cong_{\mathfrak{g}} V(\lambda),
\\
 \ch_{\gr} D(\ell, \lambda)=\ch_{\mathfrak{g}} V(\lambda)+ \sum_{\stackrel{(\mu,r)\in P^+\times \bn}{\mu< \lambda}}
\dim\Hom_{\mathfrak{g}}(V(\mu), D(\ell, \lambda)[r])q^r.\nonumber
\end{gather}
 We refer the reader to \cite{FoL} and \cite{Naoi} for a~more traditional def\/inition of a~Demazure module: namely as a~module (for a~Borel subalgebra) which is generated by an extremal vector in a~level $\ell$-highest weight module for an af\/f\/ine Kac--Moody algebra. We remark also for later use, that this approach allows one to prove that the bounds given in \eqref{demrelo} are tight, i.e.,
  \begin{gather}\label{tight}
  r<\max\{0,\lambda(h_\alpha)-d_\alpha\ell p\}\implies 0\ne (x^-_\alpha\otimes t^p)^{r+1}w_\lambda\in D(\ell,\lambda).
  \end{gather}

{\bf 2.6.} The following is now clear.
 \begin{Lemma}\label{chind}
 The elements of the set $\{{\rm{ch}}_{\gr} W_{\loc}(\lambda):\lambda\in P^+\}$ $($resp.\ $\{{\rm{ch}}_{\gr} D(\ell,\lambda): \lambda\in P^+\}$, for a~fixed $\ell\in\bn)$ are a~$\bz[q]$-basis for the subspace of $\bz[P][q]$ spanned by $\{{\rm{ch}}_{\mathfrak{g}} V(\lambda):\lambda\in P^+\}$.
 \end{Lemma}
 In particular the lemma shows that the character of a~local Weyl module is a~$\bz[q]$ linear combination of Demazure characters. We shall see in what follows that this is more than a~formal linear combination.


{\bf 2.7.}~The next result gives the relationship between $W_{\loc}(\lambda)$ and the Demazure modules $D(1,\lambda)$. It was established in \cite{CPweyl} for $\mathfrak{sl}_2$. This was then used in \cite{CL} and \cite{FoL} to prove the result for $\mathfrak{sl}_{r+1}$ and the general simply-laced case respectively. In the non-simply-laced case the result was proved in~\cite{Naoi}. A dif\/ferent approach to part (i) of this theorem along with a~generalization to the higher level Demazure modules is developed in \cite[Section~3]{CV}.
 \begin{Theorem}\label{Theorem1} Let $\lambda\in P^+$.
 \begin{enumerate}\itemsep=0pt
 \item[$(i)$] Assume that either $\mathfrak{g}$ is simply-laced or that $\lambda(h_i)=0$ for all $i\in I$ such that $\alpha_i$ is short. We have an isomorphism of $\mathfrak{g}[t]$-modules, \[W_{\loc}(\lambda)\cong D(1,\lambda).\]
 \item[$(ii)$] Let $\lambda(h_i)\ne 0$ for some $i\in I$ with $\alpha_i$ a~short root. Then we have a~decreasing sequence of graded $\mathfrak{g}[t]$-submodules \[W_{\loc}(\lambda)= W_0\supset W_1\supset\cdots\supset W_k\supset W_{k+1}=0,\] such that $W_j/W_{j+1}\cong \tau_{s_j}D(1,\mu_j)$ for some $(\mu_j,s_j)\in P^+\times\bz_+$, $0\le j\le k$.
 \end{enumerate}
\end{Theorem}

\begin{Remark} We make some comments on the proof given in \cite{Naoi} of part (ii) of the preceding theorem. An important ingredient of the proof is the following statement. Assume that $\mathfrak{g}$ is of type $\mathfrak{sl}_{n+1}$. For $\ell\ge 1$, the module $W_{\loc}(\lambda)$ (or equivalently, by part (i) of the theorem $D(1,\lambda)$), has a~f\/iltration in which the successive quotients are isomorphic to $D(\ell,\mu)$ for some $\mu\in P^+$. The result is only needed when $\ell=2$ if $n>2$ and for $\ell=2, 3$ if $n=1$. However, Naoi establishes this result for an arbitrary simply-laced Lie algebra and arbitrary $\ell$. As we remarked in the introduction the proof of this result is indirect, using a~result of A.~Joseph~\cite{Jo} for Demazure modules in the quantum case. His result uses, in a~serious way, the theory of canonical bases and is proved only for simply-laced af\/f\/ine Lie algebras, i.e., is not proved for~$A_1^1$.
 Our paper gives a~direct proof in the case of Demazure modules of af\/f\/ine $\mathfrak{sl}_2$. This is essential in Naoi's work and is also needed in the recent preprint~\cite{BMM}.
 \end{Remark}


{\bf 2.8.}~Let $V$ be a~graded $\mathfrak{g}[t]$-module. Motivated by the previous discussion, we call a~decreasing sequence \[\mathcal{F}(V)=(V= V_0\supseteq V_1\supseteq\cdots\supseteq V_k\supseteq V_{k+1}=0),\] of graded $\mathfrak{g}[t]$-submodules of $V$, a~Demazure f\/lag $\mathcal{F}(V)$ of level $\ell$, if
 \[
 V_j/V_{j+1}\cong \tau_{s_j}D(\ell,\mu_j),\qquad (\mu_j,s_j)\in P^+\times\bz_+,\qquad 0 \le j\le k.
 \]
 Let $[\mathcal{F}(V): \tau_sD(\ell,\mu)]$ be the number of times $\tau_sD(\ell,\mu)$ occurs in this f\/lag. The graded multiplicity of $D(\ell,\mu)$ in this f\/lag is def\/ined to be
 \[
 [\mathcal{F}(V) :D(\ell,\mu)]_q=\sum_{s\ge 0}[\mathcal{F}(V) : \tau_sD(\ell,\mu)]q^s.
 \] Using Lemma~\ref{chind}, one can show easily that $[\mathcal{F}(V) :D(\ell,\mu)]_q$ is independent of the choice of the f\/lag and hence from now on we will denote this polynomial by $[V:D(\ell,\mu)]_q$. The following is trivial.
\begin{Lemma}\label{elemdemflag}
Let $V$ be a~graded $\mathfrak{g}[t]$-module.
\begin{enumerate}\itemsep=0pt
\item[$(i)$] Suppose that $U$ is a~graded $\mathfrak{g}[t]$-submodule of $V$ such that $U$ and $V/U$ both admit a~Demazure flag of level $\ell$. Then $V$ also admits a~Demazure flag of level $\ell$ and \[[V: D(\ell,\mu)]_q= [U: D(\ell,\mu)]_q+ [V/U: D(\ell,\mu)]_q.\]
 \item[$(ii)$] If $V$ admits a~Demazure flag of level $\ell$, we have an equality of graded characters \[{\rm{ch}}_{\gr} V= \sum_{\mu\in P^+}[V: D(\ell,\mu)]_q{\rm{ch}}_{\gr} D(\ell,\mu).\]
 \end{enumerate}
\end{Lemma}

{\bf 2.9.}~The next result is a~reformulation of the result proved in \cite[Sections~4 and~9]{Naoi} of which the special cases will be of interest to us in this paper.
\begin{Proposition}\label{bng2} Suppose that $\mathfrak{g}$ is of type $B_n$ and that $\ell=2$ or that $\mathfrak{g}$ is of type $G_2$ and $\ell=3$. Let $\lambda\in P^+$ be such that $\lambda(h_n)=r\in\bz_+$ where $\alpha_n$ is the unique short root. Then, $W_{\loc}(\lambda)$ has a~Demazure flag of level one and for $\mu\in P^+$,
\begin{gather*}
[W_{\loc}(\lambda)): D(1, \mu)]_q=
\begin{cases}
[W^{\mathfrak{sl}_2}_{\loc}(\lambda(h_n)): D^{ \mathfrak{sl}_2}(\ell,\mu(h_n))]_q, & (\lambda-\mu)\in\bz_+\alpha_n, \\
0, & \text{otherwise}.
\end{cases}
\end{gather*}
Here $W^{\mathfrak{sl}_2}_{\loc}(\lambda(h_n))$ is the local Weyl module for $\mathfrak{sl}_2$ of weight $\lambda(h_n)$ and $D^{\mathfrak{sl}_2}(\ell,\mu(h_n))$ is defined similarly.
\end{Proposition}

We shall revisit this result later in the paper and give a~closed formula for the graded multiplicities when~$\mathfrak{g}$ is of type~$B_n$ and generating series for the numerical multiplicities ($q=1)$ when~$\mathfrak{g}$ is of type~$G_2$. We remark that the analogous result when $\mathfrak{g}$ is of type~$F_4$ or~$C_n$ is also proved in~\cite{Naoi}: here one takes $\ell=2$ and replaced~$\mathfrak{sl}_2$ by~$\mathfrak{sl}_3$ and~$\mathfrak{sl}_{n-1}$ respectively.

\section[The modules $V(\xi)$ and Demazure flags]{The modules $\boldsymbol{V(\xi)}$ and Demazure f\/lags}\label{section3}

 In this section we shall, for the most part, be interested only in the case when $\mathfrak{g}$ is $\mathfrak{sl}_2$. We shall denote the unique simple root by $\alpha$ and the corresponding fundamental weight by $\omega$. We also set $x_\alpha^+ =x$, $x_\alpha^-=y$ and $h=[x, y]$ so that $\{x,y, h\}$ is the standard basis of~$\mathfrak{sl}_2$.

{\bf 3.1.}~Recall that a~partition with $s$ parts is a~f\/inite decreasing sequence of positive integers $\xi=(\xi_1\ge\xi_2\ge \cdots\ge\xi_s>0 )$ and that $|\xi|=\sum\limits_{r\ge 1}\xi_r$. We shall sometimes also write $\xi$ as the sequence $1^{b_1}2^{b_2}\cdots s^{b_s}$, where $b_j$ is the number of times the integer $j$ occurs in~$\xi$.

 Following \cite{CV}, given a~partition $\xi$ we def\/ine a~$\mathfrak{sl}_2[t]$-module $V(\xi)$ as follows. It is the $\mathfrak{sl}_2[t]$-module generated by an element $v_\xi$ with def\/ining relations:
 \begin{gather}\label{locweylg}
 (x\otimes\bc[t])v_\xi=0,\qquad (h\otimes t^s)v_\xi=|\xi|\delta_{s,0}v_\xi,\\
 \label{integloc} (y\otimes 1)^{|\xi|+1}v_\xi=0,\\
 \label{fusrel} (x\otimes t)^s(y\otimes 1)^{s+r}v_\xi=0,
 \end{gather}
 where $s\in\bz_+$ in \eqref{locweylg} and \eqref{fusrel} holds for all $s,r\in\bn$ for which there exists $k\in\bn$ with $ s+r\ge 1+ rk+\sum\limits_{j\ge k+1}\xi_j$. It is clear that $V(\xi)$ is a~graded $\mathfrak{sl}_2[t]$-module once we set the grade of~$v_\xi$ to be zero, and that we have a~surjective map of $W_{\loc}(|\xi|\omega)\to V(\xi)$ of graded $\mathfrak{sl}_2[t]$-modules. Observe also that if~$\xi$ is the empty partition then $V(\xi)$ is the trivial $\mathfrak{sl}_2[t]$-module.

The following was proved in \cite[Proposition 2.7]{CV} and will be needed later.
\begin{Lemma} \label{defvxi}
For all $k\in\bz_+$, we have
\[
\big(y\otimes t^k\big)^{r+1}v_\xi=0,\qquad r\ge\sum_{j\ge k+1}\xi_j.
\]
\end{Lemma}

The structure of the modules~$V(\xi)$ were studied in great detail in~\cite{CV} for simple Lie algebras. In the case of~$\mathfrak{sl}_2$, these modules are known to be the fusion product of simple irreducible f\/inite-dimensional representations with dimensions~$\xi_i+1$. These were f\/irst studied in~\cite{FL} and~\cite{FF} where an alternate presentation can be found.


{\bf 3.2.}~In certain special cases the def\/ining relations of $V(\xi)$ can be greatly simplif\/ied. Thus let $\xi=1^\ell$ and, suppose that $s,r,k\in\bn$ are such that $ s+r\ge 1+ rk+\sum\limits_{j\ge k+1}\xi_j$. It is clear that $s+r\ge |\xi|+1$ and hence \eqref{fusrel} is a~consequence of \eqref{integloc}, i.e.,
\[
V(1^{\ell})\cong W_{\loc}(\ell\omega).
\]
If $\xi=(\ell)$ and if we take $s=r=1$, we f\/ind that \eqref{fusrel} holds for $k=1$, and so,
\[
(x\otimes t)(y\otimes 1)^2v_\xi=0.
\]
Using the commutation relations in $\bu(\mathfrak{sl}_2[t])$ and the def\/ining relations of $V(\xi)$ it is simple to prove that $(y\otimes t)v_\xi=0$. In turn, this implies that
\begin{gather}\label{demev}
(\mathfrak{sl}_2\otimes t\bc[t])v_\xi=0,\qquad {\rm{i.e.}}\qquad V(\xi)\cong\ev_0 V(\ell\omega),
\end{gather} as $\mathfrak {sl}_2[t]$-modules.
 More generally, the following was proved in \cite[Theorem~2]{CV}.
 \begin{Theorem}\label{demxi} Suppose that $\xi= sm^\ell $ for some $m,s\in\bn$ with $0< s\le m$ and $\ell\in\bz_+$.
 Then~$V(\xi)$ is isomorphic to the quotient of $W_{\loc}(\lambda)$ $($where $\lambda=(m\ell+s)\omega)$ obtained by imposing the additional relations
 \[
 \big(y\otimes t^{\ell+1}\big)w_\lambda=0,\qquad (y\otimes t^{\ell})^{s+1} w_\lambda =0,\qquad s<m.
 \]
 Moreover we have an isomorphism of $\mathfrak{sl}_2[t]$-modules,
 \[
 V(\xi)\cong D(m, \lambda).
 \]
\end{Theorem}

{\bf 3.3.}~Given $ m,\ell\in\bn$ with $m\ge\ell$, set \[\qbinom{m}{\ell}_q=\frac{(1-q^m)(1-q^{m-1})\cdots (1-q^{m-\ell+1})}{(1-q^{\ell})(1-q^{\ell-1})\cdots (1-q)}.\] The main result of this section is the following.
\begin{Theorem} \label{main1}\quad
\begin{enumerate}\itemsep=0pt
\item[$(i)$] For all $\xi=(\xi_1\ge\xi_2\ge\cdots\ge\xi_s>0)$ and $m\in\bn$ the module $V(\xi)$ has a~Demazure flag of level $m$ if and only if $m\ge \xi_1$. In particular, a~level $k$-Demazure module has a~Demazure flag of level $m$ if and only if $m\ge k$.
\item[$(ii)$] The local Weyl module $W_{\loc}(s\omega)$ has a~level two Demazure flag and
\[
[W_{\loc}(s\omega): D(2,r\omega)]_q= \begin{cases} q^{p\lceil s/2\rceil}\qbinom{\lfloor s/2\rfloor}{p}_q, & r=s-2p\ge 0,  \\ 0, & \text{otherwise}.\end{cases}
 \]
 \item[$(iii)$] We have an equality of characters
 \[
 \ch_{\gr}D(2, s\omega)=\sum_{p\ge 0}(-1)^{p} q^{p(\lceil{s/2}\rceil-(p-1)/2)}\qbinom{\lfloor{s/2\rfloor}}{p}_q \ch_{\gr}W_{\loc}((s-2p)\omega).\]
\end{enumerate}
\end{Theorem}

 \begin{Remark}
 Recall from equation~\eqref{demeva} that $D(m, s\omega)\cong\ev_0 V(s\omega)$ for all $m\ge s$. In particular, if we take $m=|\xi|$, a~level $m$-Demazure f\/lag is just a~Jordan--Holder series of $V(\xi)$. In the case when $\xi=1^m$ the Jordan--Holder series is known (see \cite{CL}) and that
\[
[W_{\loc}(m\omega): \ev_0 V((m-2\ell)\omega)]_q=\qbinom{m}{\ell}_q-\qbinom{m}{\ell-1}_q.
\]
 In the next section, we associate to each partition $\xi$ and an integer $m\ge\xi_1$, a~directed edge labeled graph $\mathbb H_m(\xi)$. The graded multiplicities $[V(\xi): D(m,s\omega)]_q$ can then be described by natural polynomials associated to this graph.
 The graph $\mathbb H_m(1^\ell)$ allows one to do computations with Sage to generate information on the graded multiplicities. But, even in the case when $q=1$, where one can use the OEIS database it seems dif\/f\/icult to guess the appropriate generating series for the multiplicities in general. See however Subsection~3.10, 
 where we discuss the level three Demazure f\/lag of the local Weyl module.
\end{Remark}

{\bf 3.4.}~The following is now immediate from Proposition~\ref{bng2} and Theorem~\ref{main1}.
\begin{Proposition}\label{chdem} Suppose that $\mathfrak{g}$ is of type $B_n$ and that $\lambda\in P^+$. Then we have an equality of characters \begin{gather*}
{\rm{ch_{\gr}}}W_{\loc}(\lambda)=\sum_{r\ge 0}q^{r\lceil\lambda(h_n)/2\rceil}\qbinom{\lfloor\lambda(h_n)/2\rfloor}{r}_q {\rm{ch_{\gr}}}D(1,\lambda-r\alpha_n),\\
 {\rm{ch_{\gr}}}D(1,\lambda)=\sum_{s\ge 0}(-1)^{s} q^{s(\lceil{\lambda(h_n)/2}\rceil-(s-1)/2)}\qbinom{\lfloor{\lambda(h_n)/2\rfloor}}{s}_q {\rm{ch_{\gr}}}W_{\loc}(\lambda-s\alpha_n).
 \end{gather*}
\end{Proposition}

 {\bf 3.5.}~The proof of Theorem \ref{main1} relies on a~further result proved in \cite{CV} which we now recall. Given a~partition $\xi=(\xi_1\ge\cdots\ge \xi_s>0)$ with $s$ parts we def\/ine an associated pair of partitions $\xi^\pm$ as follows. If $s=1$, then $\xi^+=\xi$ and $\xi^-$ is the empty partition. If $s>1$, then we take
 \[
 \xi^-=(\xi_1\ge\xi_2\ge \cdots \ge\xi_{s-2}\ge\xi_{s-1}-\xi_s\ge 0),
 \]
 and
 $\xi^+$ to be the unique partition associated to the $s$-tuple $(\xi_1,\dots, \xi_{s-2},\xi_{s-1}+1,\xi_{s}-1)$.
 \begin{Theorem}\label{sl2}
 Let $\xi=(\xi_1\ge\cdots\ge \xi_s>0)$ be a~partition with $s$ parts.
 For $s>1$, there exists a~short exact sequence of $\mathfrak{sl}_2[t]$-modules
 \[
 0\to\tau_{(s-1)\xi_s}V(\xi^-)\stackrel{\varphi^-}{\longrightarrow} V(\xi){\stackrel{\varphi^+}{\longrightarrow}} V(\xi^+)\to 0,
 \] and
 \[
 \varphi^+(v_{\xi})=v_{\xi^+},\qquad \varphi^-(v_{\xi^-})= \big(y\otimes t^{s-1}\big)^{\xi_s}v_\xi.
 \]
 \end{Theorem}

 {\bf 3.6.}~Def\/ine a~partial order on the set of all partitions as follows: given partitions $\xi$ and $\xi'$ with $s$ and $s'$ parts respectively, we say that $\xi>\xi'$ if and only if $s\ge s'$ and if $s=s'$, then $\xi_s>\xi'_s$ (where $\xi_s$, $\xi'_s$ is the smallest part of $\xi$ and $\xi'$ respectively). Note that $\xi>\xi^+>\xi^-$.

We f\/irst prove by an induction (given by the partial order) that if $m\ge\xi_1$, then~$V(\xi)$ has a~Demazure f\/lag of level~$m$. If $\xi=(\ell)$, then Theorem~\ref{main1}(i) is trivially true by~\eqref{demev}. Assume now that we have proved the result for all partitions with at most $(s-1)$ parts and that~$\xi$ is a~partition with~$s$ parts. The inductive hypothesis applies to~$\xi^\pm $ provided that $m\ge \xi_1^+$ and hence $\tau_{s-1}V(\xi^-)$ and $V(\xi^+)$ both have a~Demazure f\/lag of level~$m$. Lemma~\ref{elemdemflag} now shows that~$V(\xi)$ also admits a~Demazure f\/lag of level~$m$.

It remains to consider the case when $m<\xi_1^+$ which can only happen when $m=\xi_1$ and $\xi_1=\xi_{s-1}$,
 i.e., $\xi=(\xi_1^{s-1},\xi_s)$. But now, Theorem \ref{demxi} implies that $D(m, |\xi|\omega)\cong V(\xi)$ and there is nothing to prove.
 We note the following consequence of our discussion:
 \begin{gather} \label{summult}
 [V(\xi): D(m,r\omega)]_q=[V(\xi^+): D(m,r\omega)]_q+ q^{(s-1)\xi_s}[V(\xi^-): D(m,r\omega)]_q.
 \end{gather}

{\bf 3.7.}~To complete the proof of part~(i), we must show that if $V(\xi)$ has a~Demazure f\/lag of level~$m$, then $m\ge \xi_1$. This clearly follows from:
\begin{Lemma}
Suppose that $r,m\in\bz_+$ and $\xi$ is a~partition. If $D(m,r\omega)$ is a~quotient of $V(\xi)$, then $r=|\xi|$ and $m\ge\xi_1$.
\end{Lemma}

\begin{proof} Since $hv_\xi=|\xi|v_\xi$, and $D(m,r\omega)$ is a~quotient of $V(\xi)$ we must have $r=|\xi|$ and the element $v_\xi$ maps to the the generator $w_{r\omega}$ of $D(m,r\omega)$. To prove that $m\ge\xi_1$ we see from Lemma~\ref{defvxi} that
\[
(y\otimes t)^{1+\sum\limits_{j\ge 2}\xi_j}v_\xi=0,
\] and hence we must have
\[
(y\otimes t)^{1+\sum\limits_{j\ge 2}\xi_j}w_{r\omega}=0.
\] Using~\eqref{tight}, we see that this means that
\[
\sum_{j\ge 2}\xi_j\ge |\xi|-m,
\] which forces $m\ge\xi_1$ as required.
\end{proof}

{\bf 3.8.}~To prove part (ii) and for later use we establish the following result. Note that
\[
[V(\xi): D(m,r\omega)]_q\ne 0\implies (|\xi|-r)\in 2\bz_+.
\]

 \begin{Lemma}\label{mredundant}
 Given a~partition $\xi=(\xi_1\ge \xi_2\ge\cdots\ge \xi_s> 0)$ with $\xi_1=m$ and $s\ge 2$, we have \[[V(\xi): D(m, r\omega)]_q= q^{(|\xi|-r)/2}[V(\xi(m)): D(m, (r-m)\omega)]_q,\] where $\xi(m)=(\xi_2\ge\cdots\ge\xi_s> 0)$.
 \end{Lemma}

 \begin{proof} The proof is by an induction very similar to the one used in Subsection~3.6. 
 If $\xi_{s-1}=m$, then by Theorem \ref{demxi} we have
\[
 V(\xi)\cong D(m, |\xi|\omega),\qquad V(\xi(m))\cong D(m, (|\xi|-m)\omega),
\]
and there is nothing to prove. In particular induction begins at $s=2$ and we may assume also that $\xi_{s-1}<m$.
 Since
 $\xi^-= (m\ge\xi_2\ge\cdots\ge\xi_{s-1}-\xi_s)$ and $\xi^+$ is the partition associated to $(m, \xi_2,\dots, \xi_{s-1}+1, \xi_s-1)$ we get
\[
 \xi^\pm(m )=\xi(m)^\pm.
\]
Equation \eqref{summult} gives
 \begin{gather} \label{summulta}
 [V(\xi(m)): D(m,r\omega)]_q=[V(\xi^+(m)): D(m,r\omega)]_q+ q^{(s-2)\xi_s}[V(\xi^-(m)): D(m,r\omega)]_q.
 \end{gather}
Applying the induction hypothesis to $\xi^\pm$, and noting that $|\xi|=|\xi^+|$, we get
\begin{gather}\label{xim}
[V(\xi^-): D(m,r\omega)]_q= q^{(|\xi^-|-r)/2}[V(\xi^-(m)): D(m,(r-m)\omega)]_q,
\\
\label{xip}[V(\xi^+): D(m,r\omega)]_q= q^{(|\xi|-r)/2}[V(\xi^+(m): D(m,(r-m)\omega)]_q.
\end{gather}
Substituting \eqref{xim} and \eqref{xip} in \eqref{summult}, recall that $|\xi^-|=|\xi|-2\xi_s$ and then use \eqref{summulta} to complete induction.
 \end{proof}

{\bf 3.9.}~We now prove part (ii) of Theorem \ref{main1}. Since $W_{\loc}(s\omega)\cong V(1^s)$, we have
\[
\xi=1^s,\qquad \xi^+=21^{s-2},\qquad \xi^-= 1^{s-2},
\]
and hence we get
\begin{gather*}
[W_{\loc}(s\omega): D(2, r\omega)]_q = [ V(\xi^+): D(2, r\omega)]_q +q^{s-1}[W_{\loc}((s-2)\omega): D(2,r\omega)]_q\\
\qquad{}= q^ {(s-r)/2}[W_{\loc}((s-2)\omega) : D(2, (r-2)\omega)]_q + q^{s-1}[W_{\loc}((s-2)\omega): D(2,r\omega)]_q,
\end{gather*} where the second equality follows from Lemma~\ref{mredundant}. Setting
\[
p_{s,r}(q)=[W_{\loc}(s\omega): D(2, r\omega)]_q,
\] we get the following recurrence relation,
\[
p_{s,r}(q)= q^{(s-r)/2}p_{s-2,r-2}(q)+ q^{s-1}p_{s-2,r}(q),\qquad s\ge r\ge 2,
\]
with the initial conditions $p_{0,0}(q)= p_{1,1}(q)=1$ and $p_{s,r}(q)=0$ if $r>s$ or $(s-r)\notin2\bz_+$.
This recurrence has a~unique solution and a~simple calculation shows that taking
\[
p_{s, s-2\ell}=q^{\ell\lceil s/2\rceil}\qbinom{\lfloor s/2\rfloor}{\ell}_q,
\] satisf\/ies the recurrence and the initial conditions and part~(ii) is proved.
 We may regard part~(ii) as giving the change of basis matrix (see Lemma~\ref{chind})
 from the basis of local Weyl modules to the basis of level two Demazure modules. Part~(iii)
 is then just given by the inverse of this matrix and is a~straightforward calculation.

{\bf 3.10.}~We make some comments about the general case. Consider an arbitrary partition $\xi=(\xi_1\ge\cdots\ge\xi_s\ge 0)$ and suppose that we are interested in computing the graded multiplicity of $D(m,r\omega)$ for some $m\ge\xi_1$ and $r\in\bz_+$. One way to calculate this is to use Theorem~\ref{sl2} repeatedly. In the f\/irst iteration, one uses it to calculate the multiplicities of the Demazure f\/lag of level $\xi_1$, then one computes the multiplicity of $D(\xi_1+1, r\omega)$ in a~Demazure f\/lag of level $(\xi_1+1) $ in $D(\xi_1, s\omega)$ and so on. In other words, we have, that $[V(\xi): D(m,r\omega)]_q$ is equal to
\[
\sum_{\bop\in\bz_+^{m-\xi_1}}[V(\xi): D(\xi_1, p_1\omega)]_q[D(\xi_1,p_1\omega): D(\xi_1+1, p_2\omega)]_q\cdots[D(m-1: p_{m-\xi_1}\omega): D(m, r\omega)]_q.
\]
This brings us to the natural question of calculating the multiplicity occurring in a~Demazure f\/lag of $D(m-1, s\omega)$ of level $m$. In principle this can be done by using Lemma~\ref{mredundant} to set up $(m-1)$ recursive formulae. We illustrate this in the case of $m=3$ which is of interest to us and already shows the complexity of the recursions. Set
\[
p_{s,r}(q)= [D(2,s\omega): D(3,r\omega)]_q,
\] and note that $p_{s,r}=0$ if $s<r$ and also if $s$ and $r$ do not have the same parity. We shall use this freely.
Writing $s=2s_1+s_0$, $0\le s_0\le 1$ and using Lemma~\ref{mredundant} and Theorem~\ref{sl2}, we get
\begin{gather}\label{first}
p_{s,r}(q)= \begin{cases} q^{(s-r)/2}p_{s-3,r-3}(q)+ q^{(s_1+s_0-1)(2-s_0)}p_{ s-4+2s_0, r}(q), \qquad s\ge r\ge 3, \\
 q^{(s_1+s_0-1)(2-s_0)}p_{ s-4+2s_0, r}(q), \qquad s\ge 3,\quad r=0,1,2, \\
1, \qquad (r,s)\in\{(0,0), (1,1), (2,2)\},
\\ 0, \qquad {\rm{otherwise}}.\end{cases}
\end{gather}

 Setting
 \[
 P_{2r}(q,u)=\sum_{s\ge 0} p_{2s,2r}(q)u^s,\qquad P_{2r+1}(q,u)=\sum_{s\ge 0}p_{2s+1,2r+1}(q)u^{s+1},\qquad r\ge 0,
 \] we get the following equality of generating functions,
 \begin{gather*} q^{r}P_{2r}(q,u)= uqP_{2r-3}(q,qu)+u^2q^{r+2}P_{2r}\big(q, uq^2\big),\qquad r\ge 2,\\
 q^rP_{2r+1}(q,u)=u^2qP_{2r-2}(q,uq) + uq^rP_{2r+1}(q,uq),\qquad r\ge 1,
\end{gather*}
and
\begin{gather*}
P_0(q,u) = \sum_{s\ge 0}q^{2s^2}u^{2s},\qquad P_1(q,u) = \sum_{s\ge 0} q^{\frac{s(s+1)}{2}}u^{s+1},\qquad P_2(q,u) = \sum_{s\ge 0} q^{2s(s+1)}u^{2s+1}.
\end{gather*}
 If we def\/ine the partial theta series (see~\cite{AB} for instance) as
 \[
 \Theta(q,u)=\sum_{s\ge 0}q^{s^2}u^s,
 \] then note that
 \begin{gather*}
 P_0(q,u)= \Theta\big(q^2,u^2\big),\qquad P_2(q,u)= u\Theta\big(q^2,q^2u^2\big),\\
 P_1(q,u)= u\Theta\big(q^2,qu^2\big) + qu^2\Theta\big(q^2,q^3u^2\big).
 \end{gather*}
In the case when $q=1$, the recurrences can be solved and we f\/ind
\[
\big(1-u^2\big)P_{2r}= uP_{2r-3}, \qquad (1-u) P_{2r+1}= u^2P_{2r-2},
\]
or equivalently
\[
P_{r}(u)= \frac{u^3}{(1-u)(1-u^2)}P_{r-6},\qquad r\ge 6,
\] with initial conditions
\begin{gather*}P_0=\frac{1}{1-u^2},\qquad P_1=\frac{u}{1-u},\qquad P_2=\frac{u}{1-u^2},\\
 P_3=\frac{u^2}{(1-u)(1-u^2)},\qquad P_4=\frac{u^2}{(1-u)(1-u^2)},\qquad P_5=\frac{u^3}{(1-u)(1-u^2)}.
 \end{gather*}
Summarizing, we have
\begin{Proposition} The module $D(2,s\omega)$ has a~flag by level three Demazure modules and the numerical multiplicity is given by
\begin{gather*}
[D(2,s\omega): D(3,r\omega)]= \text{coefficient of} \ \ u^r\ \ \text{in} \ \
\left(\frac{u^{3r_0}P_{r_1}}{(1-u)^{r_0}(1-u^2)^{r_0}}\right),\\
 r=6r_0+r_1,\qquad 0\le r_1<6.
\end{gather*}
\end{Proposition}

Together with part (ii) of Theorem~\ref{main1} this allows us to calculate the numerical multiplicities in the f\/iltration of the local Weyl module by the level three Demazure modules and hence in turn the character of the local Weyl module when~$\mathfrak{g}$ is of type~$G_2$ in terms of Demazure modules.

{\bf 3.11.}~For arbitrary $m\in\bz_+$, one can set $p_{s,r}(q)=[D(m-1,s\omega):D(m,r\omega)]_q$ and write down recurrence relations with initial conditions which are similar to those in~\eqref{first}. To def\/ine generating functions it is not enough to just index by the weight~$s\omega$ of the level $m$-Demazure module, we have to also keep track of the smallest part of the partition, i.e., the value of $s_0$ where $s=(m-1)s_1+s_0$, $0\le s_0\le m-2$. This leads to a~large number of dependent recurrence relations and makes solving even for the numerical multiplicities very dif\/f\/icult.

{\bf 3.12.}~We conclude this section by explaining the connection of our work with Macdonald polynomials. We avoid introducing all the notation related to Macdonald polynomials, and refer the interested reader to \cite[Section 3]{CI} for a~quick exposition which is adequate for our purposes. It was shown in~\cite{S} for~$\mathfrak{sl}_r$ and in~\cite{Ion} for an arbitrary simply-laced Lie algebra that the character of the Demazure module $D(1,\lambda)$, $\lambda\in P^+$ in the level one representation of the associated untwisted af\/f\/ine Lie algebra, is given by the specialized Macdonald polynomial $P_{w_0(\lambda)}$, where $w_0$ is the longest element of the Weyl group of $\mathfrak{g}$. In the non-simply laced case however this is not necessarily true. In fact, by combining the recent work of~\cite{LNSSS} and of~\cite{Naoi} we f\/ind that the character of $W_{\loc}(\lambda)$ is given by the specialized Macdonald polynomial~$P_{w_0(\lambda)}(q)$. In other words, the character of~$W_{\loc}(\lambda)$ is always given by~$P_{w_0\lambda}(q)$. Proposition~\ref{chdem} can therefore be regarded as giving the character of the Demazure module in a~level one representation of the af\/f\/ine algebra $B_n^{(1)}$ explicitly as a~combination of Macdonald polynomials.

\section[The graph $\mathbb H_\ell(\xi)$]{The graph $\boldsymbol{\mathbb H_\ell(\xi)}$} \label{section4}

The goal of this section is to give a~combinatorial way to calculate the multiplicities of a~specif\/ic Demazure module occurring in a~Demazure f\/lag of $V(\xi)$ of level $\ell$. To do this we def\/ine a~graph~$\mathbb H_\ell(\xi)$. The vertices of the graph are a~f\/inite subset of~$\bz_+^3$ and are in bijective correspondence with the set of partitions obtained by an iteration of the operations $\xi^\pm$ (see Theorem~\ref{sl2}) with the restriction that the maximal part of the partition never exceeds~$\ell$. The edges of the graph are def\/ined to match up with the short exact sequence in Theorem~\ref{sl2} and have suitable labels. The sinks of the graph correspond to the Demazure modules of level $\ell$ which could occur with non-zero multiplicity in the f\/lag. The graded multiplicity is then obtained by summing over all ``weighted paths'' from the vertex in~$\mathbb H_\ell(\xi)$ corresponding to $\xi$ and the vertex corresponding to the desired Demazure module.

{\bf 4.1.}~Given integers $s,r\in\bn $ let $\res_r s$ be the remainder if $s$ is not divisible by $r$ and $r$ otherwise. Writing $s= (\res_rs)+ rp$ for some $p\ge 0$, we see that
\begin{gather}\label{edgeis}
s\ge r+1\implies p\ge 1\implies s\geq 2\res_rs.
\end{gather}

For $\ell\in\bn$ and a~partition $\xi=(1^{b_1}2^{b_2}\cdots )$ where $b_s=0$ if $s> \ell$ set
 \[
 \mathbb V_\ell(\xi)=\left\{(i,j,k)\in\bz_+^3: i\le\sum_{m=1}^{k-1}mb_m,\ i-jk\geq k(1-\delta_{k,\ell}),\ 2\le k\le \ell\right\}.
 \]
It is convenient to introduce the following notation:
for $\boi=(i,j,k)\in\mathbb V_\ell(\xi)$,
let $r(\boi),s(\boi)\in\bz_+$ be def\/ined by
\[
r(\boi)=\res_{k-1}(i-jk), \qquad i-jk=r(\boi)+s(\boi)(k-1).
\]

\begin{Lemma} \label{ipm}
 Let $\xi = 1^{b_1}2^{b_2}\cdots \ell^{b_{\ell}}$ and $\boi=(i,j,k)\in\mathbb V_\ell(\xi)$ be such that either $k<\ell$ or $k=\ell$ and $i-j\ell\geq\ell$. \begin{enumerate}\itemsep=0pt
 \item[$(i)$] Either $(i,j+1,k)\in\mathbb V_\ell(\xi)$ or, there exists $m(\boi)^+\in\bz_+$ with $k\le m(\boi)^+<\ell$ minimal such that $(i+b_mm, 0,m+1)\in\mathbb V_\ell(\xi)$ where $m=m(\boi)^+$.
 \item[$(ii)$]
 Either $(i-2r(\boi),j,k)\in \mathbb V_{\ell}(\xi)$ or there exists $m(\boi)^-\in\bz_+$ with $k\le m(\boi)^-<\ell$ minimal such that $(i-2r(\boi)+b_mm, 0,m+1)\in\mathbb V_\ell(\xi)$ where $m=m(\boi)^-$.
 \end{enumerate}
 \end{Lemma}

 \begin{proof}
 The proofs of part (i) and (ii) are identical. Notice that if
 \[
 i-jk\ge k(2-\delta_{k,\ell})\implies i-(j+1)k\ge k(1-\delta_{k,\ell})\implies (i,j+1,k)\in \mathbb V_{\ell}(\xi),
 \]
 resp.
 \[
 i-2r(\boi) - jk \geq k(1-\delta_{k,\ell})\implies (i-2r(\boi),j,k)\in\mathbb V_{\ell}(\xi).
 \]
Otherwise, the choice of $(i,j,k)$ forces $k<\ell$ and so $k\le i-jk<2k$.
Then \[
i+b_{\ell-1}(\ell-1)\ge 0 \qquad \text{and}\qquad i +b_{\ell-1}(\ell-1) \leq \sum_{m=1}^{\ell-1}mb_m,
\]
resp.\ by using equation \eqref{edgeis},
\begin{gather*} i-2r(\boi)+b_{\ell-1}(\ell-1) \geq jk \geq 0,
\qquad \text{and}\qquad i -2r(\boi)+b_{\ell-1}(\ell-1) \leq \sum_{m=1}^{\ell-1}mb_m),
\end{gather*}
 i.e. \[
 (i+ (\ell-1)b_{\ell-1}, 0, \ell)\in\mathbb V_{\ell}(\xi),\qquad \text{resp.} \qquad
 (i-2r(\boi)+b_{\ell-1}(\ell-1),0,\ell)\in\mathbb V_{\ell}(\xi).
 \] Hence there exists a~minimal $k\le m=m(\boi)^\pm<\ell$ such that $(i+b_mm, 0,m+1)\in\mathbb V_\ell(\xi)$ (resp.\ $(i-2r(\boi)+b_mm, 0,m+1)\in\mathbb V_\ell(\xi)$) and the lemma is established.
 \end{proof}

{\bf 4.2.}~Set \[
\boi^+=\begin{cases} (i,j+1,k), & i-jk\ge k(2-\delta_{k,\ell}), \\ (i+mb_m,0, m+1), & {\rm{otherwise}}, \quad m=m(\boi)^+
\end{cases}
\]
and
\[
\boi^-=\begin{cases} (i-2r(\boi),j,k), & i-2r(\boi)-jk\ge k(1-\delta_{k,\ell}), \\
(i-2r(\boi)+mb_m,0, m+1),& {\rm{otherwise}},\quad m=m(\boi)^-,
\end{cases}
\] where $m=m(\boi)^\pm$ is as in Lemma~\ref{ipm}. The graph $\mathbb H_\ell(\xi)$ is def\/ined to be the directed graph with vertices $\mathbb V_\ell(\xi)$ and directed edges are $\boi\to\boi^\pm$ if $\boi=(i,j,k)$ is such that either $k<\ell$ or $k=\ell$ and $i-j\ell\ge \ell$. In particular, the elements $(i,\lfloor{i/\ell}\rfloor,\ell)$ are precisely the sinks of the graph.
The label of $\boi\to\boi^+$ is zero while the label of $\boi\to\boi^-$ is
\[
\ell(\boi)= r(\boi)\left(j+b_k+\cdots+ b_{\ell}+ s(\boi) \right).
\]

{\bf 4.3.}~Let $q$ be an indeterminate. Given vertices $\boi,\boj\in\mathbb V_\ell(\xi)$ def\/ine a~polynomial $p_{\boi, \boj}(q)$ as follows. Set $p_{\boi,\boi}=1$. If $\boi$ is not a~sink, set
$p_{\boi,\boi^+}=1$, $p_{\boi,\boi^-}= q^{\ell(\boi)}$.
  Given a~directed path $\bop$ from
$(\boi=\boi_0\to\boi_1\to\cdots\to\boi_s=\boj)$ with labels $a_1,\dots ,a_s$ set $r(\bop)=\sum\limits_{m=1}^s a_m$ and def\/ine \[p_{\boi,\boj}(q)=\sum_{\bop}q^{r(\bop)},\] where the sum is over all directed paths from $\boi$ to $\boj$.
 In all other cases we take $p_{\boi,\boj}(q)=0$. The following is immediate:
 \begin{Lemma} \label{inductive}
 Given $\boi,\boj\in \mathbb V_\ell(\xi)$, we have
 \[
 p_{\boi,\boj}=p_{\boi^+,\boj}+ q^{\ell(\boi)}p_{\boi^-,\boj}.
 \]
 \end{Lemma}

 {\bf 4.4.}~We give an example. Suppose that $\ell=3$ and $\xi=1^8$. Then ${\mathbb H}_3(1^8)$ is the following edge labeled directed graph. For convenience we have color coded the edges to ref\/lect the type of edge. Thus the green (resp.\ blue) arrow is an edge corresponding to $\boi\to \boi^+$ coming from the f\/irst (resp.\ second) possibility in Lemma~\ref{ipm}(i) while the black and red arrows correspond to the two possibilities for $\boi\to\boi^-$ coming from Lemma~\ref{ipm}(ii).

\begin{figure}[t]\centering
\includegraphics[scale=0.8]{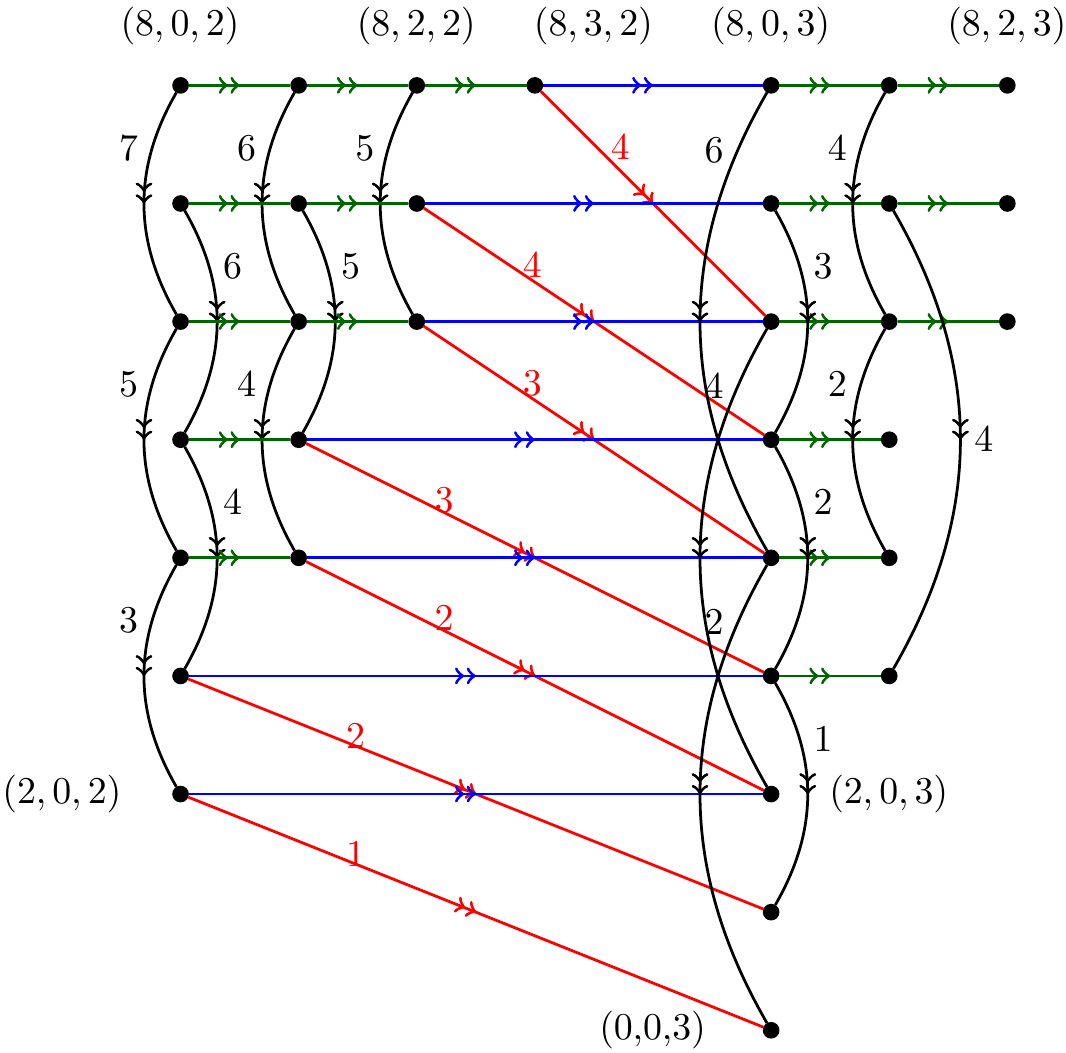}
\caption{$\mathbb H_3(1^8)$.}
\end{figure}

We give a~few examples of $p_{\boi,\boj}(q)$ for $\boi,\boj\in\mathbb V_{3}(1^8)$ which can be read of\/f easily from the graph.
\[p_{\boi,\boj}(q) =\begin{cases}
1, & \boi=(8,0,2), \quad \boj=(8,0,3),\\
q^{15}, & \boi=(8,0,2), \quad \boj=(2,0,2),\\
q^{10}+q^8, & \boi=(8,2,2), \quad \boj=(0,0,3).\\
\end{cases}\]

{\bf 4.5.}~The main result of this section is the following.
\begin{Theorem}\label{main2} Given a~partition $\xi=1^{b_1}\cdots s^{b_{s}}$ with $b_s=0$ if $s>\ell$, we have
\[
[V(\xi): D(\ell, r\omega)]_q=
\begin{cases} 0,& r<b_{\ell}\ell,\\
p_{\boi_\xi, \boj_r}(q),& \boi_\xi=(b_1,0,2), \quad \boj_r=(r-b_{\ell}\ell, \left\lfloor\frac{r}{\ell}\right\rfloor-b_{\ell},\ell).
\end{cases}
\]
\end{Theorem}

To illustrate Theorem \ref{main2} we give an example for $[V(1^8):D(3,2\omega)]_q$. In this case, $\boi_{\xi}=(8,0,2)$ and $\boj_r=(2,0,3)$. Thus,
 \begin{gather*}
[V(1^8):D(3,2\omega)]_q = p_{(8,0,2),(2,0,3)}
=q^{15}+ q^{14}+q^{13}+q^{12}+q^{11}+q^{10}+q^9+q^8= q^8[8]_q.
 \end{gather*}

The theorem is proved in the rest of the section. If $\xi=m\ell^{b_{\ell}}$ for some $0\le m<\ell$, then $V(\xi)\cong D(\ell, m+b_{\ell}\ell)$ and hence $[V(\xi): D(\ell, r)]_q=0$ if $r<b_{\ell}\ell$. Otherwise, writing $\xi=(\xi_1\ge\cdots\ge\xi_s>0)$ we see by using the def\/initions of $\xi^\pm$ that
\[
b_{\ell}=\#\{j:\xi_j=\ell\} =\#\{j:\xi^-_j=\ell\}\ \le \#\{j:\xi^+_j=\ell\}.
\]
Recall that
\[
[V(\xi): D(\ell, r\omega)]_q= [V(\xi^+): D(\ell, r\omega)]_q + q^{(s-1)\xi_s} [V(\xi^-): D(\ell, r\omega)]_q.
\]
A similar induction used in Section~\ref{section2} proves the assertion that
\[
[V(\xi): D(\ell, r\omega)]_q= 0,\qquad r<b_{\ell}\ell.
\]

{\bf 4.6.}~To prove the rest of the theorem, we relate the vertices of the graph $\mathbb H_\ell(\xi)$ to partitions and the edges to the short exact sequence in Theorem \ref{sl2}.
 Given $\boi=(i,j,k)\in\mathbb V_\ell(\xi)$, Since $i-jk\ge k(1-\delta_{k,\ell})$ we have by equation \eqref{edgeis} that $s(\boi)\ge (1-\delta_{k,\ell})$ and {\em we assume this from now on}. Let $\xi(\boi)$ be the partition
 \[
 \xi(\boi)=\begin{cases}
 r(\boi)(k-1)^{s(\boi)}k^{j+b_k}(k+1)^{b_{k+1}}\cdots\ell^{b_{\ell}}, & \boi \neq \big(i,\frac{i}{\ell},\ell\big),\\
 \ell^{b_{\ell}+\frac{i}{\ell}}, & \boi=\big(i,\frac{i}{\ell},\ell\big).
 \end{cases}
 \]

 It is trivially checked that the assignment $\boi\to \xi(\boi)$ is one-one. Note also that
 \[
 \xi(\boi_\xi)=\xi, \qquad \xi(\boj_r) =
 \begin{cases}
 r(\boi)\ell^{\left\lfloor{\frac{r}{\ell}}\right\rfloor}, & \ell \not| \, r, \\
 \ell^{\frac{r}{\ell}}, & \text{otherwise},
 \end{cases}
 \]
 where $\boi_\xi$ and $\boj_r$ are as def\/ined in Theorem \ref{main2}. In particular, \[V(\xi(\boj_r))\cong D(\ell, r\omega).\]

 \begin{Lemma}\label{lemipm}
 Let $\boi\in\mathbb V_\ell(\xi)$ and assume that either $k<\ell$ or $k=\ell$ and $i-j\ell\geq \ell$. Then,
 \[
 \xi(\boi)^\pm=\xi(\boi^\pm).
 \]
 \end{Lemma}

 \begin{proof} By def\/inition, we have
 \[
 \xi(\boi)^+=(r(\boi)-1)(k-1)^{s(\boi)-1}k^{j+b_k+1}(k+1)^{b_{k+1}}\cdots\ell^{b_{\ell}}.
 \]
 Suppose f\/irst that
 $i-jk\ge k(2-\delta_{k,\ell})$ in which case $s(\boi)\ge 2 - \delta_{k,\ell}$ and $\boi^+= (i,j+1,k)$ and we have
\begin{gather*}
 i-(j+1)k= (r(\boi)-1)+ (s(\boi)-1)(k-1).
\end{gather*}
 If $r(\boi)>1$, then it follows immediately that $i-(j+1)k$ is not divisible by $(k-1)$and hence we get
\[
r(\boi^+)=r(\boi)-1,\qquad s(\boi^+)=s(\boi)-1,
\] while if $r(\boi)=1$, then we get that
\[
r(\boi^+)=k-1,\qquad s(\boi^+)=s(\boi)-2.
\] In either case, we have now proved the assertion of the lemma.

It remains to consider the case when $k\le i-jk<2k$. We have that either $s(\boi)=2$ and $r(\boi)=1$, i.e., $i-jk=2k-1$ or $s(\boi)=1$. If $s(\boi)=2$, then $i\ge 2k-1\ge k+1$ and hence $\boi^+=(i+b_kk, 0,k+1)$. Moreover,
\[
i+b_kk= 1+ jk+ 2(k-1)+ b_kk,\qquad {\rm i.e.}, \qquad i+b_kk=k-1+ k(j+b_k+1),
\] i.e.\ $s(\boi^+)=j+b_k+1$, $r(\boi^+)=k-1$,
 and we have
\[
\xi(\boi)^+=(k-1)k^{j+b_k+1}(k+1)^{b_{k+1}}\cdots\ell^{b_{\ell}}=\xi(\boi^+).
\]
Assume now that $s(\boi)=1$. If $\res_{k-1}(i-jk)=1$, then we get
\[
\xi(\boi)^+= (\res_k k ) k^{j+b_k}(k+1)^{b_{k+1}}\cdots\ell^{b_{\ell}}.
\]
By Lemma~\ref{ipm} there exists $k\le m<\ell$ minimal such that $i+b_mm\ge (m+1)(1-\delta_{m+1,\ell})$. If $j>0$ then $m=k$, $\boi^+=(i+b_kk, 0,k+1)$ and we see that $\xi(\boi)^+= \xi(\boi^+)$ while if $j=0$ then either $k<m<\ell$ minimal with $b_m>0$ or $m=\ell$. Noting that $\res_k k=\res_m k$, it follows again that $\xi(\boi^+)=\xi(\boi^+)$.

 Finally, we must consider the case when $s(\boi)=1$ and $r(\boi)>1$. Here we have $i-jk\ge k+1$ and hence $\boi^+=(i+b_kk, 0, k+1)$. This time we just need to check that
\[
r(\boi)-1=r(\boi^+),\qquad s(\boi^+)= b_k+j+1.
\]
Since
\[
i-jk= \res_{k-1}(i-jk)+k-1\implies i+b_kk= (\res_{k-1}(i-jk)-1)+k(j+b_k+1),
\]
we see that $i+b_kk$ is not divisible by $k$ and so $0<\res_k(i+b_kk)=\res_{k-1}(i-jk)-1$ as needed.

The proof that $\xi(\boi)^-=\xi(\boi^-)$ is a~similar computation and we omit the details.
 \end{proof}

{\bf 4.7.}~Theorem \ref{main2} is clearly established if we prove the following

 \begin{Proposition}
 For all $\boi\in\mathbb V_\ell(\xi)$ and
 $\boj_r=(r-b_{\ell}\ell, \lfloor\frac{r}{\ell}\rfloor - b_{\ell}, \ell)\in\mathbb V_\ell(\xi)$, with $r\ge b_{\ell}\ell$, we have
\begin{gather*}
 [V(\xi(\boi)): V(\xi(\boj_r))]_q= p_{\boi, \boj_r}(q).
 \end{gather*}
\end{Proposition}
\begin{proof}
We proceed by induction on the length of the longest path connecting $\boi$ and $\boj_r$. If $\boi=\boj_m$ for some $m\ge b_{\ell}\ell$, then the proposition is true \[V(\xi(\boj_r))\cong D(\ell, r\omega), \qquad V(\xi(\boj_m))\cong D(\ell,m\omega).\] In other words induction begins when the length is zero. The inductive step follows by using using Theorem \ref{sl2}, Lemma~\ref{lemipm} and Lemma~\ref{inductive} which together give
\begin{gather*}
[V(\xi(\boi)): V(\xi(\boj_r)]_q= [V(\xi(\boi)^+): V(\xi(\boj_r)]_q+ q^{\ell(\boi)} [V(\xi(\boi)^-): V(\xi(\boj_r)]_q\\
\qquad {}=
[V(\xi(\boi^+)): V(\xi(\boj_r)]_q +q^{\ell(\boi)}[V(\xi(\boi^-)): V(\xi(\boj_r)]_q= p_{\boi^+,\boj_r}+ q^{\ell(\boi)}p_{\boi^-, \boj_r}= p_{\boi,\boj_r}.
\tag*{\qed}
\end{gather*}
\renewcommand{\qed}{}
\end{proof}

\subsection*{Acknowledgements}

 The authors thank S.~Viswanath for discussions regarding the graph~$\mathbb H_\ell(\xi)$ and for drawing their attention to the connection of the results of Subsection~3.10 
 to partial theta series. The f\/irst and third authors acknowledge the hospitality and excellent working conditions at the Institute of Mathematical Sciences, Chennai, India where part of this work was done. They also thank the referees of the paper for their careful reading of the paper and for their many valuable comments. The f\/irst author was partially supported by DMS-0901253 and DMS-1303052.

\pdfbookmark[1]{References}{ref}
\LastPageEnding

\end{document}